\documentclass{article}
\usepackage{amsmath}
\usepackage{latexsym}
\usepackage{amssymb}
\usepackage{graphicx}
\newtheorem{lem}{Lemma}
\newtheorem{thm}{Theorem}

\def\real{\mathbb{R}}

\def\Label#1{\label{#1}\ [\ #1\ ]\ }
\def\Label{\label}
\def\Bibitem#1{\bibitem{#1}\ [\ #1\ ]\ }
\def\Bibitem{\bibitem}
\begin{document}
\begin{center}
\bf\Large Limiting behavior of
relative R\'{e}nyi entropy\\ in a non-regular 
location shift family
\end{center}
\begin{center}
\bf\large
Masahito Hayashi\footnote{e-mail masahito@brain.riken.go.jp}\\
Laboratory for Mathematical Neuroscience, 
Brain Science Institute, RIKEN
\footnote{2-1 Hirosawa, Wako, Saitama, 351-0198, Japan}
\end{center}
\begin{abstract}
We calculate the limiting behavior of relative R\'{e}nyi entropy when the first probability distribution is close to the second one in a non-regular location-shift family which is generated by a probability distribution whose support is an interval or a half-line. This limit can be regarded as a generalization of Fisher information, and plays an important role in large deviation theory.
\end{abstract}
\section{Introduction}\Label{s1}
In a regular distribution family, Cram\'{e}r-Rao inequality holds, and the maximum likelihood estimator (MLE) converges to a normal distribution whose variance is the inverse of the Fisher information because the Fisher information converges and is well-defined in this family. However, in a non-regular location shift family which is generated by a distribution of $\real$ whose support is not $\real$ (e.g., a Weibull distribution, gamma distribution, or beta distribution), the Fisher information diverges and cannot be defined. Thus, one might think that a substitute information quantity is necessary for a discussion of the asymptotic theory. Akahira and Takeuchi \cite{AT2} proposed the limit of the Hellinger affinity $-\log \int p^{\frac{1}{2}}_{\theta}(\omega) p^{\frac{1}{2}}_{\theta+\epsilon} (\omega) \,d \omega$ as a substitute information quantity. This value is obtained by a transformation from the Hellinger distance. Moreover, Akahira \cite{Ak} proposed the relative R\'{e}nyi entropy (Chernoff's distance) $I^s(p\|q):= -\log \int p^s(\omega) q^{1-s} (\omega) \,d \omega ~(0 \,< s \,< 1)$ as a substitute information quantity for a non-regular location shift family. This quantity is linked with $\alpha$-divergence $D^{\alpha}(p\|q):=\frac{4}{1-\alpha^2} \left(1-\int_\Omega p^{\frac{1-\alpha}{2}}(\omega) q^{\frac{1-\alpha}{2}}(\omega)\,d \omega\right)$, which was introduced by Amari-Nagaoka \cite{AN} from an information geometrical viewpoint, by the monotone transformation $x \mapsto - \log \left(1- \frac{1-\alpha^2}{4}x\right)$. Since $\alpha$-divergence is a special case of
f-divergence introduced by Csisz \'{a}r\cite{Csis}, which satisfies the information processing inequality, the relative R\'{e}nyi entropy satisfies the information processing inequality
\begin{eqnarray*}
I^s ( p \| q ) \ge
I^s ( p \circ f^{-1} \| q \circ f^{-1} )
\end{eqnarray*}
for any map $f$.
Moreover, as is shown by Chernoff's formula \cite{Cher} and Hoeffding's formula \cite{Hoef}, the asymptotic error exponents in simple hypothesis testing are characterized by the relative R\'{e}nyi entropy. Thus, it can be regarded as a suitable information quantity.

As was proven by Hayashi \cite{Haya}, the upper bounds of large deviation type bounds are given by these limits of the relative R\'{e}nyi entropies. 
These upper bounds are outlined in section \ref{s4}. 
Therefore, the calculation of these limits for the non-regular location shift family is an important topic. These limits can be regarded as suitable substitutes for the Fisher information because when the Kullback-Leibler divergence is finite, the relative R\'{e}nyi entropies are connected with the Kullback-Leibler divergence by the relation
\begin{eqnarray}
D( p \| q )
=
\lim_{s \to 1} \frac{1}{s(1 - s)}
I^s( p \| q)
=
\lim_{s \to 0} \frac{1}{s(1 - s)}
I^s( q \| p) \Label{15.08} .
\end{eqnarray}
As is known, if a one-parameter distribution family ${\cal S}:=\{p_{\theta} | \theta \in \Theta \subset \real\}$ satisfies suitable regularity conditions, Kullback-Leibler divergence is closely related to the Fisher information $J_{\theta}$ defined by (\ref{h1}) as 
\begin{align}
\lim_{\epsilon \to 0}
\frac{1}{\epsilon^2} D(p_{\theta+\epsilon}\| p_{\theta})
&=
\frac{1}{2}J_{\theta}  \Label{h2} \\
J_{\theta_0}
&:=
\int_{\Omega}  
\left( \frac{\partial \log 
\frac{\partial p_{\theta}}{\partial p_{\theta_0}}(\omega)
}{\partial \theta} \right)^2
p_{\theta_0}( \,d \omega). \Label{h1}
\end{align}
However, when the support depends on the parameter $\theta$, the equation does not hold because the divergence is infinite. As was shown by Akahira \cite{Ak}, under suitable regularity conditions, the equation
\begin{align}
\lim_{\epsilon \to 0}
\frac{1}{\epsilon^2s(1-s)} I^s(p_{\theta}\| p_{\theta+\epsilon})
=
\frac{1}{2}J_{\theta}  \Label{h3}
\end{align}
holds. As the examples in sections 2 and 3 show, there are cases where relation (\ref{h3}) holds, but equation (\ref{h2}) does not. The above facts indicate that the limit of the relative R\'{e}nyi entropy is a suitable substitute for the Fisher information in a non-regular location shift family.

Moreover, in a regular family, since Fisher information is well-defined, the Riemann metric can be naturally 
defined on every tangent space. However, in a non-regular location shift family, as was pointed out by Amari \cite{Ama}, the natural metric on the tangent space is not a Riemann metric, but a general Minkowski metric. Such a manifold with a general Minkowski metric on every tangent space is called a Finsler space. Amari \cite{Ama} proposed that to treat the asymptotic behavior of the MLE, we should regard a non-regular location shift family as a Finsler space with the Minkowski metric $ F(\theta):=\lim_{\epsilon \to 0}
\frac{1}{\epsilon} 
H(p_{\theta}\|p_{\theta+\epsilon})
^{\frac{1}{\kappa}}$,
where $H$ is the Hellinger distance.
Unfortunately, the relation between the MLE and this Minkowski metric has not been adequately clarified, and the value of this Minkowski metric has not been calculated. Our result for the case $s=\frac{1}{2}$ gives the value of this Minkowski metric.

\section{Interval support case}\Label{s2}
In this section, we discuss the location shift family generated by a $C^3$ continuous probability density function $f$ whose support is an open interval $(a,b) \subset \real$. We assume conditions (\ref{con1}) and (\ref{con2}) for $f$:
\begin{align}
f_1(x):= f(a+x) &\cong  A_1 x^{\kappa_1-1} 
& \hbox{ as }x \to + 0 
\Label{con1} \\
f_2(x):= f(b-x) &\cong  A_2 x^{\kappa_2-1} 
& \hbox{ as }x \to + 0,
\Label{con2}
\end{align}
where $\kappa_1,\kappa_2 \,> 0$.
In addition, if $\kappa_i \neq 1$, we assume the following conditions:
\begin{align}
f_i'(x) &\cong A_i (\kappa_i-1) x^{\kappa_i-2}&
\hbox{ as }& x \to + 0 
\Label{con3} \\ 
f_i''(x) &\cong A_i (\kappa_i-1)(\kappa_i-2)
 x^{\kappa_i-3}& \hbox{ as }& x \to + 0 \hbox{ if }
\kappa_i \neq 2 
\Label{con4} \\
x f_i''(x) &\to 0
&\hbox{ as }& x \to +0 \hbox{ if }
\kappa_i = 2 \Label{con5}. 
\end{align}
If $\kappa_i=1$, we assume the
existence of the limits $\lim_{x \to + 0}f_i'(x)$
and $\lim_{x \to + 0} f_i''(x)$.
If $\kappa_i \,> 2$,
we assume that
\begin{align}
J_f &:=
\int_a^b f^{-1}(x) (f')^2(x) \,d x \,< 
\infty \Label{con5.1}.
\end{align}
For example, when $f$ is 
the beta distribution $f(x)=\frac{1}{B(\alpha,\beta)}
x^{\alpha-1}(1-x)^{\beta-1}$ whose support is $(0,1)$,
the above conditions are satisfied and we have
\begin{align}
\kappa_1 = \alpha, \quad \kappa_2= \beta,
\quad A_1=A_2= \frac{1}{B(\alpha,\beta)}.
\end{align}
In this paper, we denote the beta function by
$B(x,y)$.
Then, we have the following theorem.
\begin{thm}\Label{thm1}
Assume that
$\kappa:=\kappa_1 = \kappa_2$,
\begin{align}
\begin{array}{cll}
\displaystyle\lim_{\epsilon \to +0}
\frac{I^s(f_{\theta}\|f_{\theta+\epsilon})}
{\epsilon^{\kappa}}
& = \left\{
\begin{array}{l}
\displaystyle
\frac{1-\kappa}{\kappa}
(A_1 s B(s + \kappa(1-s),1-\kappa)+
A_2(1-s) B(1-s + \kappa s,1-\kappa)) \\
\displaystyle A_1s + A_2(1-s) \\
\frac{A_1 s(1-s(\kappa-1))
B(s+\kappa(1-s),2-\kappa)}{\kappa}
+\frac{A_2(1-s)(1-(1-s)(\kappa-1))
B(1-s+\kappa s,2-\kappa)}{\kappa} 
\end{array} \right.
& \begin{array}{l}
\displaystyle
0 \,< \kappa \,< 1 \\
\displaystyle \kappa = 1 \\
\displaystyle 1 \,< \kappa \,< 2 \\
\end{array} \\
\displaystyle\lim_{\epsilon \to +0}
\frac{I^s(f_{\theta}\|f_{\theta+\epsilon})}
{-\epsilon^2 \log \epsilon}
&= \displaystyle\frac{(A_1+A_2)s(1-s)}{2} 
&\displaystyle\quad\kappa =2 \\
\displaystyle\lim_{\epsilon \to +0}
\frac{I^s(f_{\theta}\|f_{\theta+\epsilon})}
{\epsilon^2}
&= \displaystyle\frac{s(1-s)}{2}J_f & \displaystyle
\quad 2 \,< \kappa ,
\end{array}  \Label{ee}
\end{align}
where $f_{\theta}(x):= f(x- \theta)$.
These convergences are uniform for $0 \,< s \,< 1$. If $\kappa_1 \,< \kappa_2$, substituting $\kappa:=\kappa_1$, $A_2:=0$, we obtain the above equations. \end{thm}
The uniformity of $0 \,< s \,<1$ is essential for the discussion in Hayashi \cite{Haya}. The case $\kappa \,> 2$ is an example where relation (\ref{h3}) holds, but relation (\ref{h2}) does not. Note that when $0 \,< \kappa \,< 2$, in general, the equation $\lim_{\epsilon \to +0} \frac{I^s(f_{\theta}\|f_{\theta+\epsilon})} {\epsilon^{\kappa}}= \lim_{\epsilon \to -0} \frac{I^s(f_{\theta}\|f_{\theta+\epsilon})} {|\epsilon|^{\kappa}}$ does not hold.

\begin{proof}
Since $I^s(f_\theta\|f_{\theta+\epsilon}) =I^s(f_{-\epsilon}\|f_0)$, Lemma \ref{thm4} yields equation (\ref{ee}). 
\end{proof}
\begin{lem}\Label{thm4}
For any $c \in (a,b)$, we define 
\begin{align*}
I^-_s(c,f,\epsilon)&:= 
\int_a^c f^{1-s}(x) f^s(x + \epsilon)
\,d x
- \int_a^c f(x) \,d x
- f(c) s \epsilon - \frac{s}{2}f'(c) \epsilon^2, \\
I^+_s(c,f,\epsilon)&:= 
\int_c^{b-\epsilon}  f^{1-s}(x) f^s(x + \epsilon)
\,d x
- \int_c^b f(x) \,d x
+ f(c) s \epsilon + \frac{s}{2}f'(c) \epsilon^2 .
\end{align*}
\begin{align}
\begin{array}{lll}
\displaystyle
\lim_{\epsilon \to +0}
\frac{I^-_s(c,f,\epsilon) }
{\epsilon^{\kappa_1}}
& =
\left\{
\begin{array}{l}
\displaystyle -\frac{1-\kappa_1}{\kappa_1}
A_1 s B(s + \kappa_1(1-s),1-\kappa_1) \\
\displaystyle -A_1s\\
\displaystyle -\frac{A_1 s(1-s(\kappa_1-1))
B(s+\kappa_1(1-s),2-\kappa_1)}{\kappa_1} 
\end{array} \right. 
&
\begin{array}{l}
\displaystyle 0 \,< \kappa_1 \,< 1 \\
\displaystyle\kappa_1 = 1 \\
\displaystyle 1 \,< \kappa_1 \,< 2 
\end{array} \\
\displaystyle\lim_{\epsilon \to +0}
\frac{I^-_s(c,f,\epsilon)}
{-\epsilon^2 \log \epsilon}
&\displaystyle = -\frac{A_1s(1-s)}{2} &\quad
\displaystyle \kappa_1 =2 \\
\displaystyle \lim_{\epsilon \to +0}
\frac{I^-_s(c,f,\epsilon)}
{\epsilon^2}
&\displaystyle = -\frac{s(1-s)}{2}J_{f,c}^- &\quad 
2 \displaystyle\,< \kappa_1 
\end{array} \Label{e1}
\end{align}
and
\begin{align}
\begin{array}{lll}
\displaystyle
\lim_{\epsilon \to +0}
\frac{I^+_s(c,f,\epsilon) }
{\epsilon^{\kappa_2}}
& = \left\{
\begin{array}{l}
\displaystyle\frac{1-\kappa_2}{\kappa_2}
-A_2(1-s) B(1-s + \kappa_2 s,1-\kappa_2))\\  
\displaystyle -A_2(1-s) \\
\displaystyle -\frac{A_2(1-s)(1-(1-s)(\kappa_2-1))
B(1-s+\kappa_2 s,2-\kappa_2)}{\kappa_2} 
\end{array} \right. &
\begin{array}{l}
\displaystyle 0 \,< \kappa_2 \,< 1 \\
\displaystyle \kappa_2 = 1 \\
\displaystyle 1 \,< \kappa_2 \,< 2 
\end{array} \\
\displaystyle \lim_{\epsilon \to +0}
\frac{I^+_s(c,f,\epsilon)}
{-\epsilon^2 \log \epsilon}
&\displaystyle = -\frac{A_2 s(1-s)}{2} 
&\displaystyle \quad \kappa_2 =2 \\
\displaystyle\lim_{\epsilon \to +0}
\frac{I^+_s(c,f,\epsilon)}
{\epsilon^2}
&\displaystyle = -\frac{s(1-s)}{2}J_{f,c}^+ 
&\displaystyle\quad 2 \,< \kappa_2 ,
\end{array}\Label{e2}
\end{align}
where
$J_f^-$ and $J_f^+$ are defined as
\begin{align*}
J_{f,c}^- :=
\int_{a}^c f^{-1}(x) (f'(x))^2 \,d x , \quad
J_{f,c}^+ :=
\int_{c}^{b} f^{-1}(x) (f'(x))^2 \,d x .
\end{align*}
These convergences are uniform for $0 \,< s \,< 1$.
\end{lem}
\section{Half-line support case}\Label{s3}
In this section, we discuss the case where the support is the half-line $(0, \infty)$ and the probability density function $f$ is $C^3$ continuous. Similarly to (\ref{con1}) and (\ref{con2}), we assume that \begin{align} f(x) \cong A x^{\kappa-1} \hbox{ as } x \to 0. \Label{con6}
\end{align}
When $\kappa \neq 1 $, we assume the following conditions: 
\begin{align}
f'(x) &\cong A_i (\kappa-1) x^{\kappa-2}&
\hbox{ as }& x \to + 0 
\Label{con7} \\ 
f''(x) &\cong A_i (\kappa-1)(\kappa-2)
 x^{\kappa-3}& \hbox{ as }& x \to + 0 \hbox{ if }
\kappa \neq 2 
\Label{con8} \\
x f''(x) &\to 0
&\hbox{ as }& x \to +0 \hbox{ if }
\kappa = 2 \Label{con9}.
\end{align}
When $\kappa=1$, we assume the existence of the limits $\lim_{x \to + 0} f'(x)$ and 
$\lim_{x \to + 0}f''(x)$.
In addition, we assume that there exist real numbers
$c \,> 0$ and $\epsilon \,> 0$ such that
\begin{align}
\int_c^\infty f^{-1}(x) (f'(x))^2 \,d x &\,< \infty 
\Label{con10}\\
\int_c^\infty 
\sup_{0\le t_1 \le \epsilon} f(x+t_1)
\sup_{0\le t_2 \le \epsilon} 
|f^{-3}(x+t_2)(f')^3 (x+t_2)| \,d x & \,< \infty 
\Label{con11}\\
\int_c^\infty 
\sup_{0\le t_1 \le \epsilon} f(x+t_1)
\sup_{0\le t_2 \le \epsilon} 
|f^{-2}(x+t_2)f'(x+t_2)f''(x+t_2)| \,d x & \,< \infty 
\Label{con12} \\
\int_c^\infty 
\sup_{0\le t_1 \le \epsilon} f(x+t_1)
\sup_{0\le t_2 \le \epsilon} 
|f^{-1}(x+t_2)f'''(x+t_2)| \,d x & \,< \infty 
\Label{con13}.
\end{align}
For example, when $f$ is Weibull distribution $f(x)= \alpha \beta x^{\alpha-1} e^{\beta x^{\alpha}}$, the above conditions are satisfied and we have 
\begin{align}
\kappa=\alpha, \quad A= \alpha \beta.
\end{align}
When $f$ is gamma distribution
$f(x)= \frac{\beta^{\alpha}}{\Gamma(\alpha)}
x^{\alpha-1}e^{\beta x}$,
the above conditions are satisfied and
\begin{align}
\kappa=\alpha, \quad A= \frac{\beta^{\alpha}}{\Gamma(\alpha)}.
\end{align}
Now, we obtain the following theorem.
\begin{thm}\Label{thm2}
We obtain
\begin{align}
\begin{array}{lll}
\displaystyle \lim_{\epsilon \to +0}
\frac{I^s(f_{\theta}\|f_{\theta+\epsilon})}
{\epsilon^{\kappa}}
& = \left\{
\begin{array}{l}
\displaystyle \frac{1-\kappa}{\kappa}
 (A s B(s + \kappa(1-s),1-\kappa) \\
\displaystyle A s \\
\displaystyle \frac{A s(1-s(\kappa-1))
B(s+\kappa(1-s),2-\kappa)}{\kappa}
\end{array} \right.
&
\begin{array}{l}
\displaystyle 0 \,< \kappa \,< 1 \\
\displaystyle \kappa = 1 \\
\displaystyle 1 \,< \kappa \,< 2 
\end{array} \\
\displaystyle \lim_{\epsilon \to +0}
\frac{I^s(f_{\theta}\|f_{\theta+\epsilon})}
{-\epsilon^2 \log \epsilon}
&\displaystyle = \frac{As(1-s)}{2} 
&\displaystyle \quad \kappa =2
\\
\displaystyle \lim_{\epsilon \to +0}
\frac{I^s(f_{\theta}\|f_{\theta+\epsilon})}
{\epsilon^2}
&\displaystyle= \frac{s(1-s)}{2}J_f 
&\displaystyle \quad 2 \,< \kappa, 
\end{array}
\Label{e3}
\end{align}
where
\begin{align}
J_f &:=
\int_0^\infty f^{-1}(x) (f')^2(x) \,d x.
\end{align}
These convergences are uniform for $0 \,< s \,< 1$.
\end{thm}
Similarly to Theorem \ref{thm1},
Theorem \ref{thm2} is proven from
Lemma \ref{thm4} and  Lemma \ref{lem2}.
\begin{lem}\Label{lem2}
For a real number $c \,> 0$ satisfying 
{\rm (\ref{con10})-(\ref{con13})},
we define
\begin{align*}
I^+_s(c,f,\epsilon)&:= 
\int_c^{\infty}  f^{1-s}(x) f^s(x + \epsilon)
\,d x
- \int_c^b f(x) \,d x
+ f(c) s \epsilon + \frac{s}{2}f'(c) \epsilon^2 .
\end{align*}
We obtain
\begin{align}
\lim_{\epsilon \to +0}
\frac{I^+_s(c,f,\epsilon)}{\epsilon^2}
=
-\frac{s(1-s)}{2}J_{f,c}^+ \Label{d1}
\end{align}
where
\begin{align*}
J_{f,c}^+
:=
\int_{c}^{\infty} f^{-1}(x) (f'(x))^2 \,d x 
\end{align*}
and the convergence of {\rm (\ref{d1})} is uniform for
$0\,<s \,<1$.
\end{lem}
\section{Relation between main results and large deviation theory}\Label{s4}
We will outline a relation between Theorems \ref{thm1} and \ref{thm2} and large deviation theory only for a location shift family $\{ f_\theta(x):= f(x-\theta)| \theta \in \real\}
$,
where $f$ satisfies the conditions given in Section \ref{s2} or Section \ref{s3}. This relation was discussed by Hayashi \cite{Haya} more precisely. As generalizations of Bahadur's large deviation type bound, we define the following quantities:
\begin{align*}
\alpha_1(\theta)&:=
\limsup_{\epsilon \to + 0}
\frac{1}{g(\epsilon)}
\sup_{\vec{T}}
\inf_{\theta-\epsilon \le \theta'\le \theta+\epsilon}
\beta(\vec{T},\theta',\epsilon) \\
\alpha_2(\theta)&:=
\sup_{\vec{T}} \liminf_{\epsilon \to + 0}
\frac{1}{g(\epsilon)}
\inf_{\theta-\epsilon \le \theta'\le \theta+\epsilon}
\beta(\vec{T},\theta',\epsilon) \\
\beta(\vec{T},\theta,\epsilon) &:=
\liminf \frac{-1}{n} \log f^n_\theta
\{ | T_n - \theta | \,> \epsilon\},
\end{align*}
where $\vec{T}=\{T_n\}$ is a sequence of estimators and $g(\epsilon)$ is chosen by
\begin{align*}
g(\epsilon)= 
\left\{
\begin{array}{ll}
\epsilon^{\kappa} & 0 \,< \kappa \,<2 \\
- \epsilon^2 \log \epsilon & \kappa = 2 \\
\epsilon^2 & \kappa \,> 2 .
\end{array}\right. 
\end{align*}
As Ibragimov and Has'minskii \cite{IH} pointed out, when KL-divergence is infinite, there exists a super efficient estimator $\vec{T}$ such that $\beta(\vec{T},\theta,\epsilon)$ and $\lim _{\epsilon \to + 0}
\frac{1}{g(\epsilon)}
\beta(\vec{T},\theta,\epsilon)$ are infinite at one point $\theta$.
Therefore, we need to take the infimum $\inf_{\theta-\epsilon \le \theta'\le \theta+\epsilon}$ into account. 
Of course, in a regular case, as was proven by Hayashi \cite{Haya}, the two bounds $\alpha_1(\theta)$ and $\alpha_2(\theta)$ coincide.

If the convergence $\lim_{\epsilon \to 0 } \frac{I^s(p_{\theta-\epsilon/2}\| p_{\theta+\epsilon/2})} {g(\epsilon)}$ is uniform for $ s \in (0,1)$ and $\theta \in K$ for any compact set $K \subset \real$, these quantities are evaluated as
\begin{align*}
\alpha_1(\theta) & \le 
\overline{\alpha}_1(\theta) := 
\left\{
\begin{array}{cl}
2^{\kappa}
\sup_{0 \,< s \,< 1}I^s_{g,\theta} & 
\hbox{ if }0\,< \kappa \,< 2\\
4
\sup_{0 \,< s \,< 1}I^s_{g,\theta} & 
\hbox{ if }\kappa\ge2 
\end{array}
\right.\\
\alpha_2(\theta)
&\le
\overline{\alpha}_2(\theta) 
:=
\left\{
\begin{array}{cl}
\sup_{0\,< s \,< 1}
\frac{I^s_{g,\theta}}{s(1-s)}
\left( s^{\frac{1}{\kappa-1}} + (1-s)^{\frac{1}{\kappa-1}} 
\right)^{\kappa-1}& \hbox{ if }
0 \,< \kappa \,< 1 \\
2 I^{\frac{1}{2}}_{g,\theta} & \hbox{ if }
\kappa = 1 \\
\inf_{0 \,< s \,< 1}
\frac{I^s_{g,\theta}}{s(1-s)}
\left( s^{\frac{1}{\kappa-1}} + (1-s)^{\frac{1}{\kappa-1}} 
\right)^{\kappa-1}& \hbox{ if }
2 \,> \kappa \,> 1 \\
\inf_{0 \,< s \,< 1}
\frac{I^s_{g,\theta}}{s(1-s)}
& \hbox{ if }2 \le \kappa ,
\end{array}
\right. 
\end{align*}
where $I^s_{g,\theta}$ are defined by
\begin{align*}
I^s_{g,\theta} &:= \lim_{\epsilon \to +0 } 
\frac{I^s(p_{\theta-\epsilon/2}\| 
p_{\theta+\epsilon/2})}
{g(\epsilon)} \quad 1 \ge s \ge 0 .
\end{align*}
Note that the uniformity of the convergence concerning $0\,< s \,<1$ is necessary for deriving the above inequalities. In Hayashi \cite{Haya}, these inequalities were proven and the attainability of bounds $\overline{\alpha}_1(\theta)$ and $\overline{\alpha}_2(\theta)$ was discussed.

\section{Conclusion}
We have calculated the limit of the relative R\'{e}nyi entropy. As mentioned in Section \ref{s4}, this calculation plays an important role in large deviation type asymptotic theory. On the other hand, we conjecture that these limits characterize the asymptotic behavior of the MLE. This relation, though, still has to be clarified.

\appendix

\section{Proof of Lemma \ref{thm4}}
\subsection{Asymptotic behavior of 
$I_s^-(c,f, \epsilon)$}\Label{s41}
In the following, when the limit $\lim_{\epsilon \to +0}g(x+\epsilon)$ ($\lim_{\epsilon \to +0}g(x-\epsilon)$) exists for a function $g$, we denote it by $g(x+0)$ ($g(x-0)$), respectively. Our situation is divided into five cases:
(i) $0\,< \kappa_1 \,<1$, (ii) $\kappa_1=1$,
(iii) $1\,< \kappa_1 \,< 2$, (iv) $\kappa_1=2$, and
(v) $\kappa_1 \,>2$.
First, we discuss cases (ii) and (v).
\begin{align}
& 
\int_a^c f^{1-s}(x) f^s(x +\epsilon) \,d x
\nonumber \\
= &
\int_{a}^c 
f^{1-s} (x) 
\left|
f^s(x+\epsilon)-
\left( f^s (x) 
+(f ^{s})'(x) \epsilon 
+(f ^{s})''(x) \frac{\epsilon^2}{2} 
\right)
\right|\,d x 
\nonumber \\
&+ \int_{a}^c \left( f (x) + f^{1-s}(x) (f ^{s})'(x) \epsilon 
+ f^{1-s}(x) (f ^{s})''(x) \frac{\epsilon^2}{2} \right)
\,d x \Label{94}
\end{align}
The first term is calculated by
\begin{align}
&\int_{a}^c \left( f (x) + f^{1-s}(x) (f ^{s})'(x) \epsilon 
+ f^{1-s}(x) (f ^{s})''(x) \frac{\epsilon^2}{2} \right)
\,d x \nonumber \\
=&
\int_{a}^c \left( f (x) + f '(x) s \epsilon 
+ 
\left( \frac{s(s-1)}{2} f^{-1}(x) (f'(x))^2 + \frac{s}{2}f''(x) \right)
\epsilon^2 \right)
\,d x \nonumber \\
=&
\int_{a}^c f (x) \,d x 
+ f(c) \epsilon s 
+ f'(c) s\frac{\epsilon^2}{2}
- f(a+0) s \epsilon 
\nonumber \\
&+ \left( 
\frac{s(s-1)}{2} \int_a^c f^{-1}(x) (f'(x))^2 \,d x
- \frac{s}{2} f'(a+0)
\right) \epsilon^2 
. \Label{941}
\end{align}
The term
\begin{align*}
\frac{1}{\epsilon^2}
\int_{a}^c 
f^{1-s} (x) 
\left|
f^s(x+\epsilon)-
\left( f^s (x) 
+(f ^{s})'(x) \epsilon 
+(f ^{s})''(x) \frac{\epsilon^2}{2} 
\right)
\right|\,d x 
\end{align*}
goes to $0$ uniformly for $0 \,< s \,< 1$ as $\epsilon \to +0$. Thus, in case (ii), since $f(a+0)=A_1$, we obtain (\ref{e1}) and the uniformity for $0\,< s \,< 1$. In case (v), since $f(a+0)=f'(a+0)=0$, we obtain (\ref{e1}) and the uniformity for $\kappa_1\,>2$.

Next, we discuss cases (i), (iii), and (iv). We can calculate $I_s^-(c,f, \epsilon)$ as  
\begin{align}
&\int_a^c f^{1-s}(x) f^s(x +\epsilon) \,d x 
\nonumber \\
=&
\int_{a+ \delta}^c 
f^{1-s} (x) 
\left|
f^s(x+\epsilon)-
\left( f^s (x) 
+(f ^{s})'(x) \epsilon 
+(f ^{s})''(x) \frac{\epsilon^2}{2} 
\right)
\right|\,d x 
\nonumber \\
&+
\int_{a+ \delta}^c 
\left( f (x) 
+ f^{1-s}(x) (f ^{s})'(x) \epsilon 
+ f^{1-s}(x) (f ^{s})''(x) \frac{\epsilon^2}{2} 
\right)\,d x 
+
\int_a^{a+\delta}
f^{1-s} (x)f ^{s}(x + \epsilon) \,d x 
\Label{1}
\end{align}
In the following, we discuss only case (i). Concerning the second term of (\ref{1}), we have
\begin{align}
& \int_{a+ \delta}^c 
\left( f (x) 
+ f^{1-s}(x) (f ^{s})'(x) \epsilon 
+ f^{1-s}(x) (f ^{s})''(x) \frac{\epsilon^2}{2} 
\right)\,d x 
+ \int_a^{a+\delta} f(x)\,d x
\nonumber \\
=&\int_{a}^c f (x) \,d x
+ \left( \int_{a+ \delta}^c f'(x) \,d x \right) s \epsilon 
+ \frac{s(s-1)}{2} \left( \int_{a+ \delta}^c f^{-1}(x) (f'(x))^2 \,d x \right)
\epsilon^2
+ \frac{s}{2} \left(\int_{a+ \delta}^c f''(x) \,d x \right)
\nonumber \\
 =&
\int_a^c f(x) \,d x  
+ (f(c)- f(a+\delta) )s \epsilon 
+ ( f'(c) - f'( a + \delta) )\frac{s}{2}\epsilon^2
+ 
\frac{s (s-1)}{2}\int_{a+\delta}^c f^{-1} (x) (f'(x))^2 \,d x
\epsilon^2  \nonumber \\
 =&
\int_a^c f(x) \,d x  
+ f(c) s \epsilon 
+ \frac{s}{2}s f'(c) \epsilon^2 \nonumber \\
& \quad 
-  f(a+ \delta) s \epsilon
- f'(a + \delta) \frac{s}{2} \epsilon^2
+ 
\frac{s(s-1)}{2} \left( \int_{a+ \delta }^c f^{-1} (x) (f'(x))^2 
\,d x \right)
\epsilon^2 
\Label{2}
\end{align}
Concerning the third term of (\ref{1}),
we can calculate
\begin{align}
&\int_a^{a+\delta}
f^{1-s} (x)f ^{s}(x + \epsilon) \,d x  - \int_a^{a+\delta} f(x)\,d x\nonumber\\
=&
\int_{a}^{a+ \delta} 
\left(f^{1-s}(x)f^s( x+\epsilon)- f(x) \right) \,d x \nonumber \\
=&
\int_{a}^{a+ \delta} \int_0^\epsilon 
f^{1-s}(x) (f^s)'(x+y) \,d y \,d x \nonumber \\
=&
\int_0^\epsilon \int_0^{\frac{\delta}{y}}
s \frac{f_1^{1-s}(yz)}{f_1^{1-s}(y(z+1))}
\frac{f_1'(y(z+1))}{f_1'(y)}\,d z
y f_1'(y) \,d y\Label{3}
\end{align}
Since 
\begin{align*}
\int_0^\infty
\frac{z^{(\kappa_1-1)(1-s)}}
{(1+z)^{(\kappa_1-1)(1-s)+2-\kappa_1}}
\,d z =
B(\kappa_1+s-\kappa_1 s,1-\kappa_1 ),
\end{align*}
using (\ref{con1}) and (\ref{con3}), we can prove that for any $\epsilon'\,>0$ real numbers $\delta \,> 0$ and $\epsilon \,>0$ exist independently for $s$ such that 
\begin{align}
\left|
\int_0^{\frac{\delta}{y}}
\frac{f_1^{1-s}(yz)}{f_1^{1-s}(y(z+1))}
\frac{f_1'(y(z+1))}{f_1'(y)}\,d z -
B(\kappa_1+s-\kappa_1 s,1-\kappa_1 ) \right|
\,< \epsilon'\Label{4}
\end{align}
for $\epsilon \,> \forall y \,> 0$.
For any $\epsilon'\,>0$, there exists a real $\epsilon \,>0$ such that
\begin{align}\Label{5}
\left|
\frac{\int_0^\epsilon y f_1'(y) \,d y}
{\epsilon^{\kappa_1}}-A_1 \frac{\kappa_1-1}{\kappa_1} \right|
\,< \epsilon' .
\end{align}
Therefore,
\begin{align}
&\left|
\frac{
\int_0^\epsilon \int_0^{\frac{\delta}{y}}
s \frac{f_1^{1-s}(yz)}{f_1^{1-s}(y(z+1))}
\frac{f_1'(y(z+1))}{f_1'(y)}\,d z
y f_1'(y) \,d y}
{\epsilon ^{\kappa_1}}
+A_1 B(\kappa_1+s-\kappa_1 s,1-\kappa_1)
\frac{s(1-\kappa_1)}{\kappa_1}
\right| \nonumber \\
 \le&
\left|
\int_0^{\frac{\delta}{y}}
s \frac{f_1^{1-s}(yz)}{f_1^{1-s}(y(z+1))}
\frac{f_1'(y(z+1))}{f_1'(y)}\,d z -
s B(\kappa_1+s-\kappa_1 s,1-\kappa_1)\right|
\frac{\int_0^\epsilon y f_1'(y) \,d y}
{\epsilon^{\kappa_1}} \nonumber \\
&\quad + 
s B(\kappa_1+s-\kappa_1 s,1-\kappa_1)
\left| \frac{\int_0^\epsilon y f_1'(y) \,d y}{\kappa_1}
+A_1 \frac{(1-\kappa_1)}{\kappa_1} \right|
{\epsilon ^{\kappa_1}}\nonumber \\
\,<& \epsilon'\left(A_1 \frac{(1-\kappa_1)}{\kappa_1}+s 
B(\kappa_1+s-\kappa_1 s,1-\kappa_1)+\epsilon'\right)
\le
\epsilon'
\left(A_1 \frac{(1-\kappa_1)}{\kappa_1}+\epsilon'
+
\sup_{0\,< s \,< 1}s B(\kappa_1+s-\kappa_1 s,1-\kappa_1)
\right).
\Label{6}
\end{align} From (\ref{1}), (\ref{2}), 
(\ref{3}), and (\ref{6}),
for any $\epsilon''\,> 0$,
there exist $\epsilon \,>0$ and $\delta \,>0$
such that
\begin{align}
& \frac{
\left|I_{s}^-(c,f,\epsilon) 
- 
B(\kappa_1+s-\kappa_1 s,1-\kappa_1)\frac{s(1-\kappa_1)}{\kappa_1}
\epsilon ^{\kappa_1}
\right|
}{\epsilon ^{\kappa_1}} \nonumber \\
\le&
\frac{1}
{\epsilon ^{\kappa_1}}
\Biggl[
\int_{a+ \delta}^c 
f^{1-s} (x) 
\left( 
f^s(x+\epsilon)-
\left( f^s (x) 
+(f ^{s})'(x) \epsilon 
+(f ^{s})''(x) \frac{\epsilon^2}{2} 
\right)
\right)\,d x  \nonumber \\
&\qquad +
\left|-  f(a+ \delta) s \epsilon
- f'(a + \delta) \frac{s}{2} \epsilon^2
+ 
\frac{s(s-1)}{2} \left( \int_{a+ \delta }^c f^{-1} (x) (f'(x))^2 
\,d x \right)
\epsilon^2 
\right| \Biggr]+\epsilon''.\Label{7}
\end{align}
The first term is less than any $\epsilon''\,> 0$ when we chose $\epsilon\,> 0$ to be sufficiently small for $\delta, \epsilon'' \,> 0$. The independence of $\epsilon \,>0$ for $0\,< s \,<1$ is shown as follows.
For any $\epsilon \,> 0$, there exists $0\le t(x,\epsilon)\le 1$ such that
\begin{align}
f^s(x+\epsilon)
-
\left( f^s(x) +(f^s)'(x) \epsilon +
(f^s)''(x) \frac{\epsilon^2}{2}\right)
=
(f^s)'''(x+ t(x,\epsilon)\epsilon
) \frac{\epsilon^3}{6}. \Label{7.1}
\end{align}
Since 
\begin{align}
(f^s)'''(x) 
=
s(s-1)(s-2)f^{s-3}(x)(f')^3(x)
+
3s(s-1) f^{s-2}(x)f'(x)f''(x)
+ s f^{s-1}(x)f'''(x),
\end{align}
we can evaluate 
\begin{align}
& \int_{a+\delta}^c
f^{1-s}(x) \left|
f^s(x+\epsilon) 
- \left( f^s(x) +(f^s)'(x) \epsilon +
(f^s)''(x) \frac{\epsilon^2}{2}\right)
\right|\,d x \nonumber \\
= &
\int_{a+\delta}^c
f^{1-s}(x) 
\left| (f^s)'''(x+ t(x,\epsilon)\epsilon) 
\frac{\epsilon^3}{6} \right| \,d x \nonumber\\
\le &
\frac{\epsilon^3}{6}
\int_{a+\delta}^c
f^{1-s}(x) f^s(x+ t(x,\epsilon)\epsilon)
\Biggl|\Biggl[
s(s-1)(s-2)f^{-3}(x+ t(x,\epsilon)\epsilon)
(f')^3(x+ t(x,\epsilon)\epsilon) \nonumber\\
&\qquad +
3s(s-1) f^{-2}(x+ t(x,\epsilon)\epsilon)
f'(x+ t(x,\epsilon)\epsilon)
f''(x+ t(x,\epsilon)\epsilon)
+ s f^{-1}(x+ t(x,\epsilon)\epsilon)
f'''(x+ t(x,\epsilon)\epsilon)\Biggr]\Biggr|
\,d x\nonumber \\
\le &
\frac{\epsilon^3}{6}
\int_{a+\delta}^c
\sup_{0\le t_1 \le \epsilon} f(x+t_1)
\Biggl[
2 \sup_{0\le t_2 \le \epsilon} 
|f^{-3}(x+t_2)(f')^3 (x+t_2)| \nonumber \\
&\qquad +
3\sup_{0\le t_3 \le \epsilon} 
|f^{-2}(x+t_3)f'(x+t_3)f''(x+t_3)| 
+
\sup_{0\le t_4 \le \epsilon} 
|f^{-1}(x+t_4)f'''(x+t_4)| \Biggr]
\,d x\Label{7.5}.
\end{align} From the $C^3$ continuity of $f$, 
the coefficient of $\epsilon^3$ at (\ref{7.5}) is finite. 
Thus, we can show the independence of $\epsilon \,>0$. We obtain (\ref{e1}) and the uniformity in case (i).

Next, we discuss cases (iii) and (iv). $\kappa_1=2$. Concerning the second term of (\ref{1}), we can calculate 
\begin{align}
&\int_{a+ \delta}^c 
\left(f (x) + f^{1-s}(x) (f ^{s})'(x) \epsilon 
+ f^{1-s}(x) (f ^{s})''(x)\frac{\epsilon^2}{2} \right)
\,d x 
+\int_a^{a+\delta} f(x)+ f^{1-s} (x)(f ^{s})'(x) \,d x \nonumber \\
& =
\int_{a}^c \left( f (x) + f'(x) s \epsilon \right)
\,d x
+ 
\frac{s(s-1)}{2} \left( \int_{a+ \delta}^c 
f^{-1}(x) (f'(x))^2 \,d x \right)\epsilon^2
+ \frac{s}{2} \left( \int_{a+ \delta}^c
 f''(x) \,d x \right)\epsilon^2
\nonumber \\
&=
\int_{a}^c f (x)\,d x + (f(c)- f(a+0)) s \epsilon 
+( f'(c) - f'(a+ \delta)) s\frac{\epsilon^2}{2}
+ \frac{s(1-s)}{2}
\left(\int_{a+\delta}^c f^{-1}(x) (f'(x))^2 \,d x\right)
\epsilon^2.
\Label{8}
\end{align}
Concerning the last term of (\ref{1}), we have
\begin{align}
&\int_a^{a+\delta}
f^{1-s} (x)f ^{s}(x + \delta) \,d x 
-\int_a^{a+\delta} f(x)+ f^{1-s} (x)(f ^{s})'(x) \,d x \nonumber \\
=&\int_{a}^{a + \delta}  \int_0^{\epsilon} \int_0^{y_1}
f^{1-s}(x) ( f^{s})''(x + y_2) \,d y_2 \,d y_1 \,d x \nonumber \\
=&
\int_0^\epsilon \int_0^{y_1} 
\int_0^{\frac{\delta}{y_2}}
\Biggl[
s \frac{f^{1-s}(y_2z)}{f^{1-s}(y_2(z+1))}
\frac{f''(y_2(z+1))}{f''(y_2)}
\frac{f''(y_2)f(y_2)}{(f')^2(y_2)}\nonumber  \\
&\qquad 
+s(s-1) 
\frac{f^{1-s}(y_2 z)}{f^{1-s}(y_2(z+1))}
\frac{f(y_2)}{f(y_2(z+1))}
\frac{(f')^2(y_2 (z+1))}{(f')^2(y_2)}
\Biggr]
\,d z 
\frac{(f')^2(y_2)}{f(y_2)} y_2 \,d y_2 \,d y_1 .\Label{9}
\end{align}

In the following, we consider only case (iii). Since
\begin{align}
\int_0^\infty
z^{(1-s)(\kappa_1-1)}
(1+z)^{s(\kappa_1-1)-2}
\,d z =
 B(1 + (1-s)(\kappa_1-1), 2-\kappa_1),
\end{align}
using (\ref{con1}), (\ref{con3}), and (\ref{con4}), we can show that for any $\epsilon' \,> 0$, there exist real numbers $\delta \,> 0$ and $\epsilon \,> 0$ such that
\begin{align}
\Biggl|\Biggl[&
\int_0^{\frac{\delta}{y_2}}
\left(
s \frac{f^{1-s}(y_2z)}{f^{1-s}(y_2(z+1))}
\frac{f''(y_2(z+1)}{f''(y_2)}
\frac{f''(y_2)f(y_2)}{(f')^2(y_2)}
+s(s-1) 
\frac{f^{1-s}(y_2z)}{f^{1-s}(y_2(z+1))}
\frac{f(y_2)}{f(y_2(z+1))}
\frac{(f')^2(y_2 (z+1))}{(f')^2(y_2)}
\right)
\,d z \nonumber \\
&- B(1 + (1-s)(\kappa_1-1), 2-\kappa_1)
\frac{s(\kappa_1-2 +(s-1)(\kappa_1-1))}{\kappa_1-1}
\Biggr]\Biggr| \,< \epsilon'\Label{10}
\end{align}
for $\epsilon \,> \forall y_2 \,> 0$.
For any $\epsilon' \,> 0$, there exists a real number $\epsilon \,> 0$ such that
\begin{align}\Label{11}
\frac{\left| \int_0^\epsilon \int_0^{y_1} 
\frac{(f')^2(y_2)}{f(y_2)}y_2 \,d y_2 \,d y_1
- \frac{\kappa_1-1}{\kappa_1}\epsilon^{\kappa_1}\right|}
{\epsilon^\kappa_1}
\,< \epsilon' .
\end{align}
Similarly to (\ref{6}), it follows from (\ref{10}) and (\ref{11}) that for any $\epsilon'' \,>0$ there exist real numbers $\delta \,> 0$ and $\epsilon\,>0$ such that
\begin{align}
\frac{
\left|
\int_{a}^{a + \delta}  \int_0^{\epsilon} \int_0^{y_1}
f^{1-s}(x) ( f^{s})''(x + y_2) \,d y_2 \,d y_1 \,d x 
+A_1 B(1 + (1-s)(\kappa_1-1), 2-\kappa_1)
\frac{s(2-\kappa_1 +(1-s)(\kappa_1-1))
\epsilon^{\kappa_1}}{\kappa_1}
\right|
}
{\epsilon^{\kappa_1}}
\,< \epsilon''.
\Label{11.1}
\end{align} From (\ref{1}), (\ref{8}), (\ref{9}), and (\ref{11}), we can evaluate
\begin{align}
&\frac{
\left|
I^-_s(c,f,\epsilon)
+ A_1 B(1 + (1-s)(\kappa_1-1), 2-\kappa_1)
\frac{s(2-\kappa_1 +(1-s)(\kappa_1-1))}{\kappa_1}
\epsilon^{\kappa_1}
\right|}
{\epsilon^{\kappa_1}} \nonumber \\
\le&
\frac{1}{\epsilon^{\kappa_1}}
\Biggl[
\int_{a+ \delta}^c 
f^{1-s} (x) 
\left|
f^s(x+\epsilon)-
\left( f^s (x) 
+(f ^{s})'(x) \epsilon 
+(f ^{s})''(x) \frac{\epsilon^2}{2} 
\right)
\right|\,d x \nonumber \\
& \qquad +
\left|
- f'(a+ \delta) s\frac{\epsilon^2}{2}
+ \frac{s(1-s)}{2}
\left(\int_{a+\delta}^c f^{-1}(x) (f'(x))^2 \,d x\right)
\epsilon^2
\right|
\Biggr]+ \epsilon''. \Label{12}
\end{align}
Note that $f(a+0)=0$. The first term is less than any $\epsilon''\,>0$ when we chose $\epsilon\,> 0$ to be sufficiently small for $\delta \,> 0$ and $\epsilon' \,> 0$. Similarly to (\ref{7}), we can show that the choice of $\epsilon \,>0$ does not depend on $0\,< s \,< 1$. Thus, we obtain (\ref{e1}) and the uniformity in case (iii).

In the following, we discuss case (iv). Using the conditions (\ref{con1}), (\ref{con3}), and (\ref{con5}), we can prove that for any $\epsilon' \,> 0$, there exist real numbers $\delta \,> 0$ and $\epsilon \,> 0$ such that
\begin{align}
&\frac{\left|
\int_0^{\frac{\delta}{y_2}}
\left(
s \frac{f^{1-s}(y_2z)}{f^{1-s}(y_2(z+1))}
\frac{f''(y_2(z+1)}{f''(y_2)}
\frac{f''(y_2)f(y_2)}{(f')^2(y_2)}
+s(s-1) 
\frac{f^{1-s}(y_2z)}{f^{1-s}(y_2(z+1))}
\frac{f(y_2)}{f(y_2(z+1))}
\frac{(f')^2(y_2 (z+1))}{(f')^2(y_2)}
\right)
\,d z 
+ s(1-s)(- \log y_2)
\right|}
{-\log y_2} \nonumber \\
& \,< \epsilon' \Label{13}
\end{align}
for $\epsilon \,> y_2 \,> 0$. For any $\epsilon' \,> 0$, there exists a real number $\epsilon \,> 0$ such that 
\begin{align}
\frac{\left| \int_0^\epsilon \int_0^{y_1} 
-\log y_2 \frac{(f')^2(y_2)}{f(y_2)}y_2 \,d y_2 \,d y_1
- A_1 (- \frac{1}{2}\epsilon^2 \log \epsilon)
\right|}
{-\epsilon^2 \log \epsilon}
\,< \epsilon'' . \Label{13.2}
\end{align}
Similarly to (\ref{6}), for any $\epsilon''\,>0$ there exist $\delta \,> 0$ and $\epsilon \,>0$ such that
\begin{align}
\frac{
\left|
\int_{a}^{a + \delta}  \int_0^{\epsilon} \int_0^{y_1}
f^{1-s}(x) ( f^{s})''(x + y_2) \,d y_2 \,d y_1 \,d x 
+ A_1 \frac{s(1-s)}{2} \epsilon^2 (- \log \epsilon)
\right|
}
{-\epsilon^2 \log \epsilon}
\,< \epsilon''
\Label{13.3}
\end{align} From (\ref{1}), (\ref{8}), 
(\ref{9}), and (\ref{13.3}),
we can evaluate this as
\begin{align}
&\frac{
\left|
I^-_s(c,f,\epsilon)
+
A_1 
\frac{s(1-s)}{2}\epsilon^2 (-\log \epsilon)
\right|}
{\epsilon^2 (-\log \epsilon)} \nonumber \\
\le
&\frac{1}{\epsilon^2 (-\log \epsilon)}
\Biggl[
\int_{a+ \delta}^c 
f^{1-s} (x) 
\left|
f^s(x+\epsilon)-
\left( f^s (x) 
+(f ^{s})'(x) \epsilon 
+(f ^{s})''(x) \frac{\epsilon^2}{2} 
\right)
\right|\,d x \nonumber \\
&\qquad +
\left|
- f'(a+ \delta) s\frac{\epsilon^2}{2}
+ \frac{s(1-s)}{2}
\left(\int_{a+\delta}^c f^{-1}(x) (f'(x))^2 \,d x\right)
\epsilon^2
\right|
\Biggr]+\epsilon''. \Label{14}
\end{align}
Note that $f(a+0)=0$.
The first term is less than any $\epsilon'' \,>0$ when we chose $\epsilon\,> 0$ to be sufficiently small for $\delta\,> 0$ and $\epsilon'' \,> 0$. Similarly to (\ref{7}), we can show that the choice of $\epsilon \,>0$ does not depend on $0\,< s \,< 1$. Thus, we obtain (\ref{e1}) and the uniformity in case (iv).

\subsection{Asymptotic behavior of $I_s^+(c,f, \epsilon)$}
As in Section \ref{s41}, our situation is divided into five cases: (i) $0\,< \kappa_2 \,<1$, (ii) $\kappa_2=1$, (iii) $1\,< \kappa_2 \,< 2$, (iv) $\kappa_2=2$,  and (v) $\kappa_2 \,>2$. First, we consider cases (ii) and (v).
\begin{align}
& \int_c^{b - \epsilon}
f^{1-s}(x)f^s(x+\epsilon) \,d x
\nonumber\\
=&
\int_c^{b - \epsilon}
f^{1-s}(x)
\left|
f^s(x+\epsilon)-
\left( f^s (x) 
+(f ^{s})'(x) \epsilon 
+(f ^{s})''(x) \frac{\epsilon^2}{2} 
\right)
\right|\,d x 
\nonumber \\
& +
\int_{c}^{b - \epsilon} \left(
f(x) +  f'(x) s \epsilon 
+ f^{-1}(x) (f')^2 (x) \frac{s(s-1)}{2}\epsilon^2
+ f''(x) \frac{s}{2} \epsilon^2
\right) \,d x \nonumber \\
=&
\int_{c}^{b}
\left( f(x) + f'(x) s \epsilon 
+ f^{-1}(x) (f')^2 (x)\frac{s(s-1)}{2}\epsilon^2
+ f''(x) \frac{s}{2} \epsilon^2 \right)
\,d x
\nonumber \\
& 
- \int_{b - \epsilon}^b
\left( f(b)  + f'(b) ( x-b) + f'(b)s \epsilon \right)
\,d x 
\nonumber \\
&+\int_c^{b - \epsilon}
f^{1-s}(x)
\left|
f^s(x+\epsilon)-
\left( f^s (x) 
+(f ^{s})'(x) \epsilon 
+(f ^{s})''(x) \frac{\epsilon^2}{2} 
\right)
\right|\,d x 
\nonumber \\
& 
-\int_{b - \epsilon}^b
\left(
f(x) +  f'(x) s \epsilon 
+ f^{-1}(x) (f')^2 (x) \frac{s(s-1)}{2}\epsilon^2
+ f''(x) \frac{s}{2} \epsilon^2
\right)
-\left( f(b-0)  + f'(b-0) ( x-b) + f'(b-0)s 
\epsilon \right)
 \,d x. \nonumber 
\end{align}
The first and second terms are calculated as
\begin{align}
&\int_{c}^{b}
\left( f(x) + f'(x) s \epsilon 
+ f^{-1}(x) (f')^2 (x)\frac{s(s-1)}{2}\epsilon^2
+ f''(x) \frac{s}{2} \epsilon^2 \right)
\,d x
\nonumber \\
&
\quad - \int_{b - \epsilon}^b
\left( f(b-0)  + f'(b) ( x-b) + f'(b-0)s \epsilon \right)
\,d x \nonumber \\
=&
\int_{c}^{b}
f(x) \,d x
+ (f(b-0) - f(c) )s \epsilon 
+ (f'(b-0) - f'(c) )\frac{s}{2}\epsilon^2
\nonumber \\
& \quad 
+ \left(\int_{c}^{b}  f^{-1}(x) (f')^2 (x)\,d x \right)
\frac{s(s-1)}{2}\epsilon^2
- \int_{b - \epsilon}^b
\left( f(b-0)  + f'(b-0) ( x-b) + f'(b)s \epsilon \right)
\,d x \nonumber \\
=&
\int_{c}^{b}
f(x) \,d x  - f(c) s \epsilon
- f'(c) \frac{s}{2}\epsilon^2
+ f(b) s \epsilon 
+ f'(b) \frac{s}{2}\epsilon^2
\nonumber \\
&\quad 
+ \left(\int_{c}^{b}  f^{-1}(x) (f')^2 (x)\,d x \right)
\frac{s(s-1)}{2}\epsilon^2
- f(b-0)  \epsilon 
+ f'(b-0) \frac{\epsilon^2}{2}
- f'(b-0)s \epsilon^2 \nonumber \\
=&
\int_{c}^{b}
f(x) \,d x  - f(c) s \epsilon
- f'(c) \frac{s}{2}\epsilon^2
\nonumber \\
&\quad 
+ f(b-0) (s-1) \epsilon 
+ f'(b-0) (1-s) \frac{\epsilon^2}{2}
+ \left(\int_{c}^{b}  f^{-1}(x) (f')^2 (x)\,d x \right)
\frac{s(s-1)}{2}\epsilon^2.
\nonumber
\end{align}
The term
\begin{align*}
&\frac{1}{\epsilon^2}\Biggl[
\int_c^{b - \epsilon}
f^{1-s}(x)
\left|
f^s(x+\epsilon)-
\left( f^s (x) 
+(f ^{s})'(x) \epsilon 
+(f ^{s})''(x) \frac{\epsilon^2}{2} 
\right)
\right|\,d x \\
& \qquad
+\int_{b - \epsilon}^b
\left| \left(
f(x) +  f'(x) s \epsilon 
+ f^{-1}(x) (f')^2 (x) \frac{s(s-1)}{2}\epsilon^2
+ f''(x) \frac{s}{2} \epsilon^2
\right)
-\left( f(b)  + f'(b) ( x-b) + f'(b)s \epsilon \right)
\right| \,d x \Biggr]
\end{align*}
goes to $0$ uniformly for $0 \,< s \,< 1$ as $\epsilon \to +0$. In case (ii), the $C^3$ continuity of $f$ and the existence of $f_2'(0)$ and $f_2''(0)$ guarantee 
\begin{align*}
\int_a^b f^{-1}(x) (f')^2(x) \,d x \,< 
\infty .
\end{align*}
Thus, from the existence of $f_2'(0)$ and the relation $f(b-0)=A_2$, we obtain (\ref{e2}) and the uniformity in case (ii). From (\ref{con5.1}) and the relations $f(b-0)=f'(b-0)=0$, we obtain (\ref{e2}) and the uniformity in case (v).

Next, we consider cases (i), (iii), and (iv). 
\begin{align}
& \int_c^{b - \epsilon}
f^{1-s}(x)f^s(x+\epsilon) \,d x
\nonumber\\
=&
\int_c^{b - \delta}
\left|
f^s(x+\epsilon)-
\left( f^s (x) 
+(f ^{s})'(x) \epsilon 
+(f ^{s})''(x) \frac{\epsilon^2}{2} 
\right)
\right|\,d x 
\nonumber \\
& +
\int_c^{b - \delta}f(x) + \epsilon f^{1-s}(x) (f^s)'(x) 
+ \frac{\epsilon^2}{2} f^{1-s}(x) (f^s)''(x) \,d x
\nonumber \\
& 
+ \int_{b - \delta}^{b- \epsilon} f^{1-s}(x) f^s(x+ \epsilon) \,d x 
. \Label{15}
\end{align}
In the following, we discuss only case (i).
\begin{align}
&\int_c^{b - \delta}f(x) + \epsilon f^{1-s}(x) (f^s)'(x) 
+ \frac{\epsilon^2}{2} f^{1-s}(x) (f^s)''(x) \,d x
+\int_{b-\delta+\epsilon}^b f(x)\,d x
\nonumber \\
=&
\int_c^b f(x) \,d x
- \int_{b- \delta}^{b- \delta+\epsilon}
f( x) \,d x \nonumber \\
& + 
\left( \int_c ^{b- \delta} f^{1-s} (x)
(f^s)'(x) \,d x \right) \epsilon
+ 
\left( \int_c^{b- \delta}
f^{1-s}(x) (f^s)''(x) \,d x \right)
\frac{\epsilon^2}{2} \nonumber \\
= &
\int_c^b f(x) \,d x
- \int_{b- \delta}^{b- \delta+\epsilon}
f( x) \,d x \nonumber
+ \left( \int_c ^{b- \delta} f'(x) \,d x \right)s \epsilon
\nonumber \\
& + 
\left( \frac{s(s-1)}{2}
\int_c^{b- \delta}
(s-1) f^{-1}(x)(f'(x))^2 \,d x 
+ \frac{s}{2} \int_c^{b- \delta}f''(x) \,d x \right)
\epsilon^2 \nonumber \\
=&
\int_c^b f(x) \,d x
- \int_{b- \delta}^{b- \delta+\epsilon}
f( x) \,d x + ( f( b - \delta) - f(c)) s\epsilon 
+ 
(f'(b - \delta) - f'(c) )
\frac{s}{2}\epsilon^2 
\nonumber \\
& +
\frac{s( s-1)}{2}
\left( \int_c^{b- \delta}
f^{-1}(x)(f'(x))^2 \,d x \right)
\epsilon^2 .
\Label{16}
\end{align}
Letting $z:= \frac{b-x}{y}$, we have
\begin{align}
& \int_{b - \delta}^{b- \epsilon} f^{1-s}(x) f^s(x+ \epsilon) \,d x 
-\int_{b-\delta+\epsilon}^b f(x)\,d x \nonumber\\
=& \int_{b - \delta + \epsilon }^b 
\left( f^{1-s}(x - \epsilon) f^s(x) - f(x) \right) \,d x 
\nonumber\\
=&\int_{b - \delta + \epsilon }^b 
f^s(x) 
\int_{- \epsilon}^0 (f^{1-s})'(x+y) \,d y \,d x  
\nonumber\\
=&
\int_0^\epsilon\int_0^{\frac{\delta-\epsilon}{y}}
\frac{f_2^s(y z)}{f_2^s(y(z+1))}
\frac{f_2'(y(z+1))}{f_2'(y)}
\,d z y f_2'(y) \,d y.\Label{17}
\end{align}
Similarly to (\ref{6}), we can prove that for any $\epsilon'' \,> 0$, there exist real numbers $\delta \,> 0$ and $\epsilon \,> 0$ such that
\begin{align}
\left|
\frac{
\int_0^{\frac{\delta-\epsilon}{y}}
\frac{f_2^s(y z)}{f_2^s(y(z+1))}
\frac{f_2'(y(z+1))}{f_2'(y)}
\,d z
}
{\epsilon ^{\kappa_2}}
+B(\kappa_2+s-\kappa_2 s,1-\kappa_2)
\frac{s(1-\kappa_2)}{\kappa_2}
\right| \,< \epsilon''
\Label{18}
\end{align}
for $\epsilon \,> \forall y \,> 0$. Therefore, from (\ref{15}), (\ref{16}), (\ref{17}), and (\ref{18}), similarly to (\ref{7}) we can prove that for any $\epsilon' \,> 0$,
an real number $\epsilon \,> 0$ exists 
independently for $s$, such that 
\begin{align}
& \frac{
\left|I_{s}^+(c,f,\epsilon) 
+ A_2
B(\kappa_2+(1-s)-\kappa_2 (1-s),1-\kappa_2)\frac{
(1-s)(1-\kappa_2)}{\kappa_2}
\epsilon ^{\kappa_2}
\right|
}{\epsilon ^{\kappa_2}} 
\,< \epsilon'' . \Label{19}
\end{align}
Thus, we obtain (\ref{e2}) and the uniformity in case (i).

Next, we consider cases (iii) and (iv). Concerning the second term of (\ref{15}), we have
\begin{align}
&
\int_c^{b - \delta}f(x) + \epsilon f^{1-s}(x) (f^s)'(x) 
+ \frac{\epsilon^2}{2} f^{1-s}(x) (f^s)''(x) \,d x
+ \int_{b-\delta + \epsilon}^b  f(x) - f^s(x) (f^{1-s})'(x) \epsilon \,d x 
\nonumber \\
=&
\int_c^b \left( f(x) + f'(x) s \epsilon \right)\,d x
+ \left(\int_c^{b-\delta} f^{1-s}(x) (f^s)''(x) \,d x\right)
\frac{\epsilon^2}{2}
\nonumber \\
& \quad 
- \int_{b - \delta}^b \left( f(x) + f'(x) s \epsilon \right)\,d x 
+ \int_{b-\delta + \epsilon}^b 
\left( f(x) - f^s(x) (f^{1-s})'(x)\epsilon \right)
\,d x \nonumber\\
=&
\int_c^b f(x) \,d x
+ f(b-0)s \epsilon - f(c)s \epsilon 
+ \left(\int_c^{b-\delta} f^{1-s}(x) (f^s)''(x) \,d x\right)
\frac{\epsilon^2}{2}
\nonumber \\
& \quad 
- \int_{b - \delta}^b \left( f(x) + f'(x) s \epsilon \right)\,d x 
+ \int_{b-\delta + \epsilon}^b 
\left( f(x) - f^s(x) (f^{1-s})'(x)\epsilon \right)
\,d x 
. \Label{20}
\end{align}
We can evaluate this as
\begin{align}
& \left| - \int_{b - \delta}^b \left( f(x) + f'(x) s \epsilon \right)\,d x 
+ \int_{b-\delta + \epsilon}^b 
\left( f(x) - f^s(x) (f^{1-s})'(x)\epsilon \right)
\,d x \right|\nonumber\\
=&\left|
- \int_{b-\delta}^{b-\delta+\epsilon} f(x) \,d x
- \int_{b-\delta}^{b-\delta+\epsilon} f'(x)s \epsilon \,d x
- \int_{b-\delta+\epsilon}^b
f'(x) \epsilon \,d x  \right|\nonumber\\
=&
\left|\int_{b-\delta}^{b-\delta+\epsilon}
f(b-\delta+\epsilon) -f(x)
-f'(x)s \epsilon
\,d x \right| \nonumber\\
\le &
\int_{b-\delta}^{b-\delta+\epsilon}
\left|
f(b-\delta+\epsilon) -f(x)\right|+
|f'(x)|s \epsilon
\,d x \nonumber\\
\le &
\max_{0 \le t \le1}
|f'(b-\delta+\epsilon t) |\frac{3}{2} \epsilon^2 
. \Label{20.1}
\end{align}
Concerning the third term of (\ref{15}), we have
\begin{align}
&\int_{b - \delta}^{b- \epsilon} f^{1-s}(x) f^s(x+ \epsilon) \,d x  
+ \int_{b-\delta + \epsilon}^b - f(x) + f^s(x) (f^{1-s})'(x) \epsilon \,d x 
\nonumber \\
=&\int_{b-\delta + \epsilon}^b f^s(x) \left(
f^{1-s} (x- \epsilon) - f^{1-s} (x) 
+ (f^{1-s})'(x) \epsilon \right) \,d x \nonumber \\
=&
\int_{b-\delta + \epsilon}^b f^s(x) 
\left(- \int_{- \epsilon}^0 (f^{1-s})'(x + y_1) - 
(f^{1-s})'(x) \,d y_1
\right) \,d x \nonumber \\
=&
 \int_{b-\delta + \epsilon}^b f^s(x) 
\left(- \int_{- \epsilon}^0 \int_0^{y_1} 
(f^{1-s})''(x + y_2) \,d y_2 \,d y_1
\right) \,d x \nonumber \\
=&
\int_0^\epsilon\int_0^{y_1}\int_0^{\frac{\delta-\epsilon}
{y_2}}
\Biggl[
(1-s) \frac{f^{s}_2(y_2z)}{f^{s}_2(y_2(z+1))}
\frac{f''_2(y_2(z+1))}{f''(y_2)}
\frac{f''_2(y_2)f(y_2)}{(f'_2)^2(y_2)} \nonumber \\
&\quad +s(s-1) 
\frac{f_2^{s}(y_2 z)}{f_2^{s}(y_2(z+1))}
\frac{f_2(y_2)}{f_2(y_2(z+1))}
\frac{(f'_2)^2(y_2 (z+1))}{(f'_2)^2(y_2)}
\Biggr]
\,d z 
\frac{(f'_2)^2(y_2)}{f_2(y_2)} y_2 \,d y_2 \,d y_1 
\Label{21}
\end{align}
Similarly to (\ref{11.1}), in case (iii), we can prove that for any $\epsilon'' \,> 0$ there exist real numbers $\delta\,>0$ and $\epsilon \,> 0$ such that
\begin{align}
\frac{
\left|
 \int_{b-\delta + \epsilon}^b f^s(x) \left(- \int_{- \epsilon}^0 \int_0^{y_1} (f^{1-s})''(x + y_2) \,d y_2 \,d y_1 \right) \,d x +A_2 B(1 + s(\kappa_2-1), 2-\kappa_2) \frac{(1-s)(2-\kappa_2 +s(\kappa_2-1)) \epsilon^{\kappa_2}}{\kappa_2}
\right|
}
{\epsilon^{\kappa_2}}
\,< \epsilon''
\Label{22}
\end{align}
Similarly to (\ref{12}), from (\ref{15}), (\ref{20}), (\ref{20.1}), (\ref{21}), and (\ref{22}), we can prove that for any $\epsilon'' \,> 0$ there exists a real number $\epsilon \,>0$ such that
\begin{align}
\frac{
\left|
I^+_s(c,f,\epsilon)
+
A_2 B(1 + s(\kappa_2-1), 2-\kappa_2)
\frac{(1-s)(2-\kappa_2 +s(\kappa_2-1))}{\kappa_2}
\epsilon^{\kappa_2}
\right|}
{\epsilon^{\kappa_2}} 
\,< \epsilon''. \Label{23}
\end{align}
Thus, we obtain (\ref{e2}) and the uniformity i 
in case (iii). Similarly to (\ref{13.3}), in case (iv), we can prove that
\begin{align}
\frac{
\left|
 \int_{b-\delta + \epsilon}^b f^s(x) 
\left(- \int_{- \epsilon}^0 \int_0^{y_1} 
(f^{1-s})''(x + y_2) \,d y_2 \,d y_1
\right) \,d x 
+ A_2 \frac{s(1-s)}{2} \epsilon^2 (- \log \epsilon)
\right|
}
{\epsilon^2(-\log \epsilon)}
\,< \epsilon''.
\Label{24}
\end{align}
Similarly to (\ref{14}), from (\ref{15}), (\ref{20}), (\ref{20.1}), (\ref{21}), and (\ref{24}), we have
\begin{align}
\frac{
\left|
I^+_s(c,f,\epsilon)
+
A_2
\frac{s(1-s)}{2}\epsilon^2 (-\log \epsilon)
\right|}
{\epsilon^2 (-\log \epsilon)} 
\,< \epsilon''. \Label{25}
\end{align}
Thus, we obtain (\ref{e2}) and the uniformity in case (iv).

\section{Proof of Lemma \ref{lem2}}
We can calculate
\begin{align}
&\int_c^\infty
f^{1-s}(x)f^s(x+\epsilon) \,d x \nonumber\\
=&
\int_c^\infty
f^{1-s}(x) \left(
f^s(x+\epsilon) 
- \left( f^s(x) +(f^s)'(x) \epsilon +
(f^s)''(x) \frac{\epsilon^2}{2}\right)
\right)\,d x \nonumber\\
&\qquad +\int_c^\infty
f^{1-s}(x) \left( f^s(x) +(f^s)'(x) \epsilon +
(f^s)''(x) \frac{\epsilon^2}{2}
\right)\,d x  . \Label{51}
\end{align}
The second term of (\ref{51}) is calculated as
\begin{align}
\int_c^\infty
& f^{1-s}(x) \left( f^s(x) +(f^s)'(x) \epsilon +
(f^s)''(x) \frac{\epsilon^2}{2}
\right)\,d x\nonumber\\
=&
\int_c^\infty
f(x) +f'(x)s \epsilon 
+ f^{-1}(x)(f')^2(x)\frac{s(s-1)}{2}
\epsilon^2
+f''(x)\frac{s}{2}\epsilon^2 \,d x  \nonumber\\
=&
\int_c^\infty f(x) \,d x 
+\int_c^\infty f^{-1}(x)(f')^2(x)\frac{s(s-1)}{2}
\,d x
- f(c)s \epsilon 
-f'(c)\frac{s}{2}\epsilon^2 \,d x  \Label{52}
\end{align}
Similarly to (\ref{7.5}), we can evaluate the first term of (\ref{51}) by
\begin{align}
& \int_c^\infty
f^{1-s}(x) \left|
f^s(x+\epsilon) 
- \left( f^s(x) +(f^s)'(x) \epsilon +
(f^s)''(x) \frac{\epsilon^2}{2}\right)
\right|\,d x \nonumber \\
\le &
\frac{\epsilon^3}{6}
\int_c^\infty
\sup_{0\le t_1 \le \epsilon} f(x+t_1)
\Biggl[
2 \sup_{0\le t_2 \le \epsilon} 
|f^{-3}(x+t_2)(f')^3 (x+t_2)| \nonumber \\
&\qquad +
3\sup_{0\le t_3 \le \epsilon} 
|f^{-2}(x+t_3)f'(x+t_3)f''(x+t_3)| 
+
\sup_{0\le t_4 \le \epsilon} 
|f^{-1}(x+t_4)f'''(x+t_4)| \Biggr]
\,d x\Label{53}.
\end{align}
Conditions (\ref{con11}) - (\ref{con13}) guarantee that the coefficient of (\ref{53}) is finite. From (\ref{51}), (\ref{52}), and (\ref{53}), we obtain (\ref{d1}) and the uniformity for $0 \,< s \,< 1$.

\end{document}